\documentclass[english]{article}
\usepackage[T1]{fontenc}
\usepackage[latin9]{inputenc}
\usepackage{babel}
\usepackage{amsmath}
\usepackage{amssymb}
\usepackage[unicode=true]
 {hyperref}

\makeatletter
\newcommand{\lyxaddress}[1]{
	\par {\raggedright #1
	\vspace{1.4em}
	\noindent\par}
}

\makeatother

\begin{document}
\title{A Tauberian characterization of the Riemann hypothesis through the
floor function}
\author{Benoit Cloitre}
\maketitle

\lyxaddress{\begin{center}
benoit.cloitre@proton.me
\par\end{center}}
\begin{abstract}
We introduce a new Tauberian framework through the theory of \textquotedbl regular
arithmetic functions\textquotedbl . This allows us to establish a
characterization of the Riemann hypothesis by linking the floor function
to the distribution of nontrivial zeros of the Riemann zeta function.
We thereby obtain a novel Tauberian equivalence of the Riemann hypothesis,
extending classical Tauberian theorems beyond their traditional confinement
to the prime number theorem. We further uncover connections to combinatorial
number theory and set the groundwork for a \textquotedbl combinatorial
Tauberian theory\textquotedbl , highlighting the broader applicability
of regular arithmetic functions.
\end{abstract}

\paragraph*{Keywords}

Number theory, Tauberian theory, Analytical number theory, Combinatorial
number theory, Regular arithmetic functions, Functions of good variation,
Riemann hypothesis, Floor function, Zeta function, Möbius function,
Hardy-Littlewood-Ramanujan criterion, Arithmetic asymptotic formulas,
Ingham summation method, Dirichlet series, Elementary methods.

\paragraph*{Mathematics Subject Classification (2020)}

\paragraph*{Primary: }

40E05, 11M26, 11B75 

\paragraph*{Secondary: }

11A05, 11A25, 11N37, 05A19, 11N05

\section*{Introduction}

The understanding of our main result requires familiarity with the
concepts of regular arithmetic functions. To achieve this, we begin
by introducing the essential definitions that form the theoretical
backbone of our framework, ensuring both clarity and self-containment.
While our work addresses deep questions such as the Riemann hypothesis,
the definitions themselves remain accessible to any reader with a
standard undergraduate background in analysis.\\

The theory of regular arithmetic functions (RAF) represents a framework
for studying functions of two variables, $G:\mathbb{N^{\ast}}\times\mathbb{N^{\ast}}\longmapsto\mathbb{R}$,
that satisfy specific properties depending on the value of a real
parameter $\beta$. These properties are tied to an unknown sequence
$\left(a_{n}\right)_{n\geq1},$ satisfying $a_{1}$=1, which is uniquely
determined through the exact relation for $n\geq2$

\begin{equation}
\sum_{k=1}^{n}a_{k}G\left(n,k\right)=n^{-\beta}.
\end{equation}
This framework constitutes a distinct branch of Tauberian theory \cite{key-20},
providing asymptotic formulas for the partial sums $\sum_{n\leq x}a_{n}$
based on the notion of regularity index of the RAF. \\

Although difference equations have been extensively studied, as documented
in \cite{key-1}, and equation (1) shares similarities with the well-researched
Volterra equation \cite{key-10}, our approach uncovers nontrivial
connections between the arithmetic and analytic properties of $G$.
As we will demonstrate, this framework sheds new light on longstanding
problems, including the Riemann hypothesis and other challenging questions.
\\

We formalize the definitions underpinning the RAF theory. These concepts
are closely interrelated: the arithmetic Mellin transform serves as
the foundation for defining regular arithmetic functions, which, in
turn, include functions of good variation as a special case. We begin
with the fundamental notion of arithmetic Mellin transform: 

\paragraph*{Arithmetic Mellin transform}

We say that $G$ has an arithmetic Mellin transform if, for $\Re z<0$,
the following limit exists:
\[
G^{\star}(z):=\lim_{n\rightarrow\infty}\left(-\frac{z}{n}\right)\sum_{k=1}^{n}G\left(n,k\right)\left(\frac{k}{n}\right)^{-z-1}.
\]

The arithmetic Mellin transform of $G$, denoted $G^{\star}$, is
the analytic continuation of this limit to the complex plane, possibly
excluding a countable set of poles or singularities, and possibly
possessing a meromorphic border. $G^{\star}$ must admit a removable
singularity at the origin.

\paragraph*{Regular arithmetic function}

We say that $G:\mathbb{N^{\ast}}\times\mathbb{N^{\ast}}\longmapsto\mathbb{R}$
is a regular arithmetic function (RAF) with regularity index $\alpha(G)\in\mathbb{R}$
if it has an arithmetic Mellin transform and, $(1)$ being satisfied,
if it exhibits the following asymptotic properties depending on the
value of $\beta\in\mathbb{R}$ relative to $\alpha(G)$:\\

\subparagraph*{Case 1: $\beta<\alpha(G)$}

In this case, 
\[
\sum_{n\le x}a_{n}\sim\frac{1}{G^{\star}(\beta)}x^{-\beta}\quad(x\to\infty).
\]

\subparagraph*{Case 2: $\beta\ge\alpha(G)$}

In this case, 
\[
\sum_{n\le x}a_{n}=\mathcal{O}(x^{-\alpha(G)+\varepsilon})\quad(x\to\infty).
\]

\paragraph*{Function of good variation}

A function $g$ defined on $]0,1]$ is said to be a function of good
variation (FGV) of regularity index $\alpha(g)$, if the associated
function $G$, defined by

\[
G(n,k)=g\left(\frac{k}{n}\right),
\]
is a RAF with index $\alpha\left(G\right)=\alpha(g)$ as per the definition
above. Additionally, we denote $g^{\star}=G^{\star}$ , the arithmetic
Mellin transform of $G$, and we will simply refer to $g^{\star}$
as the Mellin transform of $g$. It is generally given by the classical
integral:
\[
g^{\star}(z)=\int_{0}^{1}g(t)t^{-z-1}dt.
\]

At the heart of our approach lies the Ingham function $\Phi(x):=x\left\lfloor x^{-1}\right\rfloor $defined
on $(0,1]$ which appeared implicitly in Ingham's Tauberian theorem
(\cite{key-17}), which states that for a sequence $a_{n}$ satisfying
the one-sided condition $na_{n}\geq-C$, we have

\[
\sum_{k=1}^{n}a_{k}\Phi\left(\frac{k}{n}\right)=1+o(1)\Rightarrow\sum_{k=1}^{n}a_{k}\rightarrow1\ \left(n\rightarrow\infty\right).
\]
This is the case $\beta=0$ of equation (1), when $G(n,k)=\Phi\left(\frac{k}{n}\right)$,
with an error term. Like the Wiener-Ikehara theorem \cite{key-16},
Ingham's Tauberian theorem yields the prime number theorem (see also
\cite{key-20} p.$110$ for a proof).\\

However, despite multiple investigations \cite{key-18,key-25,key-27,key-28,key-34},
including earlier work by Wintner \cite{key-33} and attempts from
Segal to connect it to the Riemann hypothesis \cite{key-29}, Ingham's
summation method has not yet yielded significant breakthroughs in
approaching the Riemann hypothesis. Importantly, it has not been employed
in the manner we propose here. \\

In Tauberian theory, the Ingham function, when composed with the reciprocal
map, connects also to what some authors call the Pólya kernel $k(t)=\left\lfloor t\right\rfloor t^{-1}$
defined on $[1,\infty)$ (see for instance \cite{key-4}). Unfortunately,
the classical Mellin transform of the Pólya kernel

\[
\ensuremath{\mathcal{M}[k](s)=\frac{\zeta(1+s)}{1+s}}
\]

restricts applications to $\Re s\geq0$, where $\zeta(1+s)$ doesn't
vanish, making it difficult to enter the critical strip.\\

Studying the Ingham function on $(0,1]$, however, reveals its full
potential. Its Mellin transform $\Phi^{\star}$ (as defined above),
relates to $\zeta(1-s)$ and thereby to the zeros of $\zeta(s)$ within
the critical strip, when $\Re s\geq0$. This deceptively simple shift
in Tauberian theory brings a paradigm-changing insight, enabling progress
on problems ranging from the Riemann hypothesis to our resolution
of the open problem in \cite{key-19}, which will be detailed in a
forthcoming paper, with the promise of further applications yet to
be explored.\\

Having introduced the key concepts of RAF theory and highlighted the
potential of studying the Ingham function on $(0,1]$, we can now
present our main contribution: a Tauberian equivalence of the Riemann
hypothesis that connects arithmetic properties of the floor function
with the distribution of non-trivial zeros of the zeta function. This
equivalence appears to be new, as it is not referenced, directly or
indirectly, in \cite{key-6,key-7}.

\subsubsection*{Main result}

The following statements are equivalent: 
\begin{description}
\item [{(i)}] The Riemann hypothesis is true. 
\item [{(ii)}] The Ingham function is a FGV with regularity index $\alpha\left(\Phi\right)=1/2$.\\
\end{description}
The structure of the paper is as follows:

\paragraph*{Section 1}

We establish the existence of RAFs/FGVs to legitimate our Riemann
hypothesis equivalence approach. The section presents illustrative
examples of both FGV and RAF functions that are not FGV, thus demonstrating
the richness of these concepts. 

\paragraph*{Section 2}

We explore some properties of $\Phi$. In particular, we will show
that its Mellin transform is given by $\Phi^{\star}(z)=-\frac{z}{1-z}\zeta\left(1-z\right)$,
connecting the Ingham function to the Riemann zeta function. We will
also show that $\Phi$ is connected to the fundamental theorem of
arithmetic for the case $\beta>0$, which, as we shall see, corresponds
to the non-trivial case of our equivalence.

\paragraph*{Section 3}

We establish some lemmas needed for the proof of our main result.
In particular we give estimates for $\sum_{n\leq x}a_{n}$ when $\beta\leq0$
using elementary arithmetical techniques. When $\beta>0$ we give
estimates for partial sums of the generalised Jordan functions $J_{-\beta}(n)=\sum_{d\mid n}d^{-\beta}\mu\left(\frac{n}{d}\right)$,
requiring the use of Perron's formula. 

\paragraph*{Section 4}

We prove the main result by a careful case-by-case analysis of the
values of $\beta$. We will also see in passing that the sole knowledge
that $\Phi$ is a FGV with a strictly positive regularity index $\left(\alpha\left(\Phi\right)>0\right)$
allows us to recover the prime number theorem.

\paragraph*{Section 5}

We explore connections between RAF theory and combinatorial number
theory through generalizations of equation (1) with slowly varying
functions. Under the Riemann hypothesis, these generalizations can
provide asymptotic formulas for sums like $\sum_{n\leq x}\mu\left(6n\right)$
based on combinatorial properties of the Möbius function and the floor
function, giving birth to combinatorial Tauberian theory.

\paragraph*{Section 6}

Beyond the equivalence between the Riemann hypothesis and the regularity
property of the Ingham function, we develop a broader framework which
is outlined in this section. Through generalizations of equation (1)
and the arithmetic Mellin transform, we introduce the concepts of
stable and balanced RAFs with respect to another function. This new
perspective suggests our equivalence might open a promising path toward
a better understanding of the Riemann hypothesis in the Tauberian
setting.\\

Before delving into these developments, it is essential to understand
how our approach fundamentally differs from previous Tauberian methods
in approaching the Riemann hypothesis. This analysis will help clarify
both the novelty and the potential of our framework.

\paragraph*{On the limits of previous Tauberian methods in approaching the Riemann
hypothesis}

Classical Tauberian theorems, such as those by Ingham or using the
Polyà kernel described above, along with their generalizations, are
deeply rooted in arithmetic, maintaining strong connections through
Dirichlet series and arithmetical functions like the Möbius or the
von Mangoldt functions. However, their analytic reach is confined
to the region $\Re z\geq1$ where the zeta function does not vanish.
This limitation prevents them from directly accessing the critical
strip and they cannot prove anything beyond the prime number theorem.\\

On the other hand, Tauberian approaches such as the Nyman-Beurling
criterion \cite{key-23} and other closure problems (see \cite{key-7},
chapter 3, Banach and Hilbert spaces methods, pp.23-36 and chapter
8 with Salem's equivalence 8.4 and Levinson's equivalence 8.5) are
predominantly analytic, with their core directly tied to the behavior
of the zeta function in the critical strip. These approaches can be
considered Tauberian, although this is less apparent than with theorems
like those of Ingham or Ikehara, which are closer to Tauber's original
theorem. Indeed, these closure problems are connected to another Tauberian
theorem by Wiener (\cite{key-32}, pp. 1-100, \cite{key-20}, chapter
II, pp. 65-115). Yet, these methods lack any significant arithmetic
properties, which to us weakens their applicability in number-theoretic
contexts. As Levinson remarked in \cite{key-21}:\\

\textit{``The trivial character of all such sufficiency proofs seems
to indicate that,if the Riemann hypothesis is true, the closure theorems
do not seem to be a very promising direction to pursue}''.\\

This assessment remains certainly relevant today, as no significant
progress has been made in these approaches, despite decades of research.
While there is some persistent optimism in the field and some attempts
at arithmetizing the Nyman-Beurling criterion to make sufficiency
proofs more challenging (see Baez-Duarte result, \cite{key-7}, page
36), this arithmetization doesn't really have a connection with multiplicative
number theory, which for us represents the arithmetic core of the
Riemann hypothesis.\\

In contrast, our equivalence maintains both the arithmetical and analytical
aspects of the Riemann hypothesis. It effectively bridges these two
previously disconnected Tauberian approaches, opening a third, more
promising path. Indeed, while both classical Tauberian theorems and
the Nyman-Beurling approach have their distinct limitations, our method
synthesizes their strengths: it combines the arithmetic depth of Ingham's
Tauberian theorem with the analytic potential of Nyman-Beurling theorem
to probe the critical strip. 

While the connection of our approach to Ingham's Tauberian theorem,
from which it originates, should be evident to the reader, its relationship
to Nyman's theorem might be less obvious. However, this becomes clear
when comparing the Mellin transforms at the heart of each approach.
The Nyman-Beurling criterion relies crucially on the classical Mellin
transform:

\[
\int_{0}^{\infty}\rho\left(\frac{1}{t}\right)t^{z-1}dt=-\frac{1}{z}\zeta(z),\ 0<\Re z<1,
\]
where $\rho\left(.\right)$ is the fractional part function, whereas
our approach leverages on this form of the Mellin transform:

\[
\int_{0}^{1}\Phi\left(t\right)t^{-z-1}dt=-\frac{z}{1-z}\zeta(1-z),\ \Re z<0,
\]
where $\Phi\left(.\right)$ is the Ingham function. These Mellin transforms
are intimately connected to the zeta function and, through the functional
equation of the zeta function, they are essentially equivalent regarding
the location of nontrivial zeros in the critical strip, thanks to
analytic continuation for the latter. Section 6 will also show how
RAF theory is part of functional analysis by extending the study to
function spaces, reinforcing the analogy with Nyman's criterion.

\section{Examples of regular arithmetic functions}

We present here several instances of regular arithmetic functions
(RAF). We begin with functions of good variation (FGV), which form
a particular class of RAF expressible as $G(n,k)=g(k/n)$, before
examining RAF that do not admit such a representation. These examples
demonstrate concretely that such functions exist, with a complete
proof given for our first example, while proofs for subsequent examples
will be developed elsewhere. 

The reader may observe that the proof provided here is quite involved
despite the apparent simplicity of the example. This complexity hints
at the depth of RAF theory, where even elementary functions can exhibit
asymptotic behavior determined by subtle interactions between arithmetic
and analysis.

\subsection{Smooth and bounded FGV}

For $\lambda\in]0,1[$, consider the affine function $g:[0,1]\rightarrow\mathbb{R}$
defined by

\[
g(x)=\left(1-\lambda\right)x+\lambda.
\]
This function has the Mellin transform given by

\[
g^{\star}(z)=\lambda-\left(1-\lambda\right)\frac{z}{1-z},
\]
and $g$ is a FGV with index $\alpha(g)=\lambda,$ which corresponds
to the only zero of $g^{\star}$. 

\subsubsection*{Proof}

To prove this constructively, we will need the following lemma, which
can be immediately proved by induction. 

\paragraph*{Lemma}

If $\left(x_{n},u_{n},v_{n}\right)_{n\geq1}$ are sequences satisfying
for $n\geq2$ the recurrence relation $x_{n}=u_{n}x_{n-1}+v_{n}$,
we have for $n\geq1$ 
\[
x_{n}=\left(\prod_{k=1}^{n}u_{k}\right)\left(\frac{x_{1}-v_{1}}{u_{1}}+\sum_{k=1}^{n}\frac{v_{k}}{\prod_{i=1}^{k}u_{i}}\right).
\]

Then we get

\[
\sum_{k=1}^{n}a_{k}g\left(\frac{k}{n}\right)=n^{-\beta}\Rightarrow\left(1-\lambda\right)\sum_{k=1}^{n}ka_{k}+\lambda n\sum_{k=1}^{n}a_{k}=n^{1-\beta}.
\]

Therefore, denoting $A(n)=\sum_{k=1}^{n}a_{k}$ we have $A(1)=a_{1}=1$
and for $n\geq2$ 

\[
A(n)=\left(1-\frac{\lambda}{n}\right)A(n-1)+B(n),
\]
where $B(n)=\frac{1}{n}\left(n^{1-\beta}-(n-1){}^{1-\beta}\right)$.
The lemma then yields
\begin{equation}
A(n)=\left(\prod_{k=1}^{n}1-\frac{\lambda}{k}\right)\left(\sum_{k=1}^{n}\frac{B(k)}{\prod_{i=1}^{k}1-\frac{\lambda}{i}}\right).
\end{equation}
Since $0<\lambda<1$, the value $\Gamma\left(1-\lambda\right)$ is
well defined and non-zero, so that the Euler-Gauss formula for the
Gamma function \cite{key-12}:

\[
\Gamma(s)=\lim_{n\to\infty}\frac{n!}{s(s+1)(s+2)\cdots(s+n)}n^{s},
\]
allows us to have the following two asymptotic formulas

\begin{equation}
\prod_{k=1}^{n}\left(1-\frac{\lambda}{k}\right)\sim\frac{n^{-\lambda}}{\Gamma\left(1-\lambda\right)}\ \left(n\rightarrow\infty\right),
\end{equation}

\begin{equation}
\frac{B(k)}{\prod_{i=1}^{k}1-\frac{\lambda}{i}}\sim\left(1-\beta\right)\Gamma\left(1-\lambda\right)k^{\lambda-\beta-1}\ \left(k\rightarrow\infty\right),
\end{equation}
By summing the equivalent (4), we then obtain, depending on the values
of $\beta$

\begin{equation}
\beta<\lambda\Rightarrow\sum_{k=1}^{n}\frac{B(k)}{\prod_{i=1}^{k}1-\frac{\lambda}{i}}\sim\frac{\left(1-\beta\right)\Gamma\left(1-\lambda\right)}{\lambda-\beta}n^{\lambda-\beta}\ \left(n\rightarrow\infty\right),
\end{equation}

\begin{equation}
\beta=\lambda\Rightarrow\sum_{k=1}^{n}\frac{B(k)}{\prod_{i=1}^{k}1-\frac{\lambda}{i}}\sim\left(1-\beta\right)\Gamma\left(1-\lambda\right)\log n\ \left(n\rightarrow\infty\right),
\end{equation}

\begin{equation}
\beta>\lambda\Rightarrow\lim_{n\rightarrow\infty}\sum_{k=1}^{n}\frac{B(k)}{\prod_{i=1}^{k}1-\frac{\lambda}{i}}=\ell_{\beta}<\infty.
\end{equation}
Next the Mellin transform of $g$ has the following formula for $\Re z<0$
\[
g^{\star}(z)=-z\int_{0}^{1}\left(\left(1-\lambda\right)t+\lambda\right)t^{-z-1}dt=\frac{z}{z-1}\left(1-\lambda\right)+\lambda
\]
Therefore, by injecting formulas (3) and (5)(6)(7) into (2) according
to the cases, we obtain the following asymptotic formulas for $A(n)$
as $n$ tends to infinity

\[
\beta<\lambda\Rightarrow A(n)\sim\frac{\left(1-\beta\right)}{\left(\lambda-\beta\right)}n^{-\beta}=\frac{1}{g^{\star}\left(\beta\right)}n^{-\beta},
\]

\[
\beta=\lambda\Rightarrow A(n)\sim\left(1-\lambda\right)n^{-\lambda}\log n,
\]

\[
\beta>\lambda\Rightarrow A(n)\sim\frac{\ell_{\beta}}{\Gamma\left(1-\lambda\right)}n^{-\lambda}.
\]

Which gives in summary

\[
\beta<\lambda\Rightarrow A(n)\sim\frac{1}{g^{\star}\left(\beta\right)}n^{-\beta}\ \left(n\rightarrow\infty\right)
\]

\[
\beta\geq\lambda\Rightarrow A(n)\ll n^{-\lambda+\varepsilon}\ \left(n\rightarrow\infty\right)
\]
And $g$ is a FGV with regularity index $\alpha(g)=\lambda$ according
to the definition of a FGV given in the introduction.

\subsection{FGV not bounded at zero}

With $\lambda\in]0,1]$ define $g:]0,1]\rightarrow\mathbb{R}$ as

\[
g(x)=1-\lambda\log x
\]
which has the Mellin transform

\[
g^{\star}(z)=\frac{\lambda-z}{z^{2}}.
\]
Then $g$ is an FGV of regularity index $\alpha(g)=\lambda,$ which
corresponds to the only zero of $g^{\star}$.

\subsection{Discontinuous FGV}

For an integer $\lambda\geq2$, define the function $g:]0,1]\rightarrow\mathbb{R}$
as follows

\[
g(x)=x\lambda^{\left\lfloor -\frac{\log x}{\log\lambda}\right\rfloor }.
\]
Similar to the Ingham function, which has an infinite number of discontinuities
at $x=k^{-1},\ k\geq2$, this function exhibits discontinuities at
$x=2^{-k},\ k\geq1,$ and has the Mellin transform given by

\[
g^{\star}(z)=\frac{z}{z-1}\left(\frac{\lambda^{z-1}-1}{\lambda^{z}-1}\right),
\]
and $g$ is a FGV with index $\alpha(g)=1,$ which corresponds to
the smallest real part of the zeros of $g^{\star}$, all of which
are located on the line $x=1$.

\subsection{RAF not FGV which is a rational function}

Let $x,y$ be two distinct positive real numbers, and define the function
$G$ by 

\[
G(n,k)=\frac{n+k+x}{n+k+y}.
\]
The arithmetic Mellin transform is simply given by $G^{\star}(z)=1$
and $G$ is a RAF with index $\alpha\left(G\right)=2$. 

\subsection{Some remarks on the regularity index}

In FGV examples 1.1, 1.2, and 1.3, we observe that the regularity
index is given by the smallest real part of the zeros of the Mellin
transform. In example 1.3, if $\lambda>1$ is not integer, it turns
out that $g$ is still a FGV but the value of the index is not always
equal to $1$. We can thus show that if $\lambda=\sqrt{2}$ then $g$
is an FGV with index $\frac{1}{2}$ which is quite surprising at first
glance. Similarly, in example 1.4 of a RAF that is not a FGV, the
index cannot be given by the Mellin transform, which is constant.\\

This highlights the regularity index as a central and intrinsic characteristic
of RAF, independent of the zeros of the Mellin transform in certain
cases. Its role as a bridge between the arithmetic and analytic properties
of RAF underscores its fundamental importance in this theory. It can
be viewed as a mathematical object that, sometime, connects the arithmetic
properties of a RAF $G$ with the analytic properties of its Mellin
transform $G^{\star}$. As we will see, this remarkable connection
emerges strikingly for the Ingham function, where the regularity index
unveils profound arithmetic-analytic relationships between the floor
function and the zeta function.

\section{Properties of the Ingham function}

This section presents two results: first, we give a formula for the
Mellin transform of the Ingham function, then we show that the Ingham
function satisfies a criterion connecting asymptotic behaviour of
the sequence $a_{n}$ to the fundamental theorem of arithmetic.

\subsection{The Mellin transform of the Ingham function}

The function $\Phi(x)=x\left\lfloor \frac{1}{x}\right\rfloor $ possesses
a Mellin transform $\Phi^{\star}$ that connects to the Riemann zeta
function. Indeed, for $\Re z<0$, we have 
\[
\Phi^{\star}\left(z\right)=\frac{z}{z-1}\zeta\left(1-z\right)
\]
which analytically continues to $\mathbb{C}$ except for $z=1$, where
it has a simple pole.

\subsubsection*{Proof}

For $\Re z<0$, by applying the dominated convergence theorem, and
using $\Phi\left(t\right)t^{\Re\left(-z-1\right)}$ as an integrable
dominating function, we obtain 

\[
\Phi^{\star}\left(z\right):=\lim_{n\rightarrow\infty}\left(\frac{-z}{n}\right)\sum_{k=1}^{n}\Phi\left(\frac{k}{n}\right)\left(\frac{k}{n}\right)^{-z-1}=-z\int_{0}^{1}\Phi\left(t\right)t^{-z-1}dt.
\]
Next we get

\[
\int_{0}^{1}\Phi\left(t\right)t^{-z-1}dt=_{t\rightarrow1/x}\ \int_{1}^{\infty}\Phi\left(1/x\right)x^{z-1}dx=\sum_{n=1}^{\infty}n\int_{n}^{n+1}x^{z-2}dx
\]

\[
=\frac{1}{z-1}\sum_{n=1}^{\infty}n\left(\frac{1}{(n+1)^{1-z}}-\frac{1}{n^{1-z}}\right)
\]

\[
=\frac{1}{z-1}\sum_{n=1}^{\infty}\left(-\frac{1}{(n+1)^{1-z}}+\frac{n+1}{(n+1)^{1-z}}-\frac{n}{n^{1-z}}\right)
\]

\[
=\frac{-1}{z-1}\sum_{n=0}^{\infty}\frac{1}{(n+1)^{1-z}}=\frac{-1}{z-1}\zeta\left(1-z\right),
\]
and the result follows.$\square$\\

Note that this formula naturally yields $\zeta\left(1-z\right)$ rather
than $\zeta\left(1+z\right)$ or $\zeta\left(z\right)$. This perspective
shift highlights the potential of Ingham\textquoteright s summation
method to transcend the prime number theorem, creating a pathway to
the Riemann hypothesis for $\beta>0$ in equation (1), where $G(n,k)=\Phi\left(\frac{k}{n}\right)$,
as we will see. \\

Having examined an analytic property of the Ingham function, we now
turn to a more arithmetic one.

\subsection{The Hardy-Littlewood-Ramanujan criterion}

While the definitions of RAF and FGV (refer to the Introduction) do
not explicitly impose conditions on the sequence $a_{n}$, the interaction
with recurrence relation (1) implicitly enforces specific behavior
on $a_{n}$. This hidden behavior, mirroring the Tauberian conditions
in classical Tauberian theorems (such as in Ingham Tauberian theorem
where the condition $na_{n}\geq-C$ is needed), is exemplified by
the Hardy-Littlewood-Ramanujan (HLR) criterion. The HLR criterion
draws its name from its connection to the renowned Tauberian conditions
established by Hardy and Littlewood in their seminal work \cite{key-13},
as well as Ramanujan's condition on the size of coefficients for Dirichlet
series in the context of the Selberg class \cite{key-30}. \\

\paragraph*{Definition $2.2$ }

A FGV $g$ is said to satisfy the HLR criterion if we have:

\[
\forall\beta\geq0,\ \sum_{k=1}^{n}a_{k}g\left(\frac{k}{n}\right)=n^{-\beta}\Rightarrow\forall\varepsilon>0,\ na_{n}=O\left(n^{\varepsilon}\right).
\]

\subsection{HLR criterion for the Ingham function}

We begin by establishing a lemma concerning Möbius convolution for
multiplicative functions - a result that will crystallize how the
asymptotic behavior of $a_{n}$ deeply reflects the fundamental theorem
of arithmetic.

\subsubsection*{Lemma 2.3}

If $u$ is multiplicative and satisfies $0<u_{n}\leq1$ for $n\geq1$,
then

\[
\sum_{d\mid n}\mu\left(\frac{n}{d}\right)u_{d}=O(1).
\]

\paragraph*{Proof of lemma 2.3}

Define $b(n):=\sum_{d\mid n}\mu\left(\frac{n}{d}\right)u_{d}$ which
is multiplicative since $u$ is multiplicative. This multiplicative
function satisfies

\[
b\left(p^{v}\right)=u_{p^{v}}-u_{p^{v-1}}
\]
at prime powers. Next, using the prime factorization $n=\prod p_{i}^{\alpha_{i}}$
and exploiting the bound 
\[
-1\leq u_{p_{i}^{\alpha_{i}}}-u_{p_{i}^{\alpha_{i}-1}}\leq1,
\]
we obtain

\[
\sum_{d\mid n}\mu\left(\frac{n}{d}\right)u_{d}=\prod\left(u_{p_{i}^{\alpha_{i}}}-u_{p_{i}^{\alpha_{i}-1}}\right)=O(1).
\]

$\square$

\paragraph*{Theorem $2.3$ }

The function $\Phi$ satisfies the HLR criterion, specifically

\paragraph*{
\[
\forall\beta\protect\geq0,\ \sum_{k=1}^{n}a_{k}\Phi\left(\frac{k}{n}\right)=n^{-\beta}\Rightarrow na_{n}=O\left(1\right).
\]
}

\subparagraph*{Remark }

The HLR criterion, expressed in the form $na_{n}=O(1)$, stands as
a notable attribute of Ingham's function since the asymptotic behavior
of $na_{n}$ in Ingham's Tauberian theorem has been the subject of
dedicated research, as evidenced by the following result. Consider
the equation

\[
\sum_{k=1}^{n}a_{k}\Phi\left(\frac{k}{n}\right)=\ell+o(1),
\]
then one can show that $na_{n}=o\left(\log\log n\right)$ is the best
possible asymptotic result \cite{key-11}. Therefore, equation $(1)$
without error term or a version with a suitably small error term is
necessary for $na_{n}=O(1)$.

\subsubsection{Proof of theorem $2.3$}

The proposition holds for $\beta=0$ because in this case, $a_{1}=1$
and $a_{n}=0$ for $n\geq2$. Now let $\beta>0$ and suppose

\[
\sum_{k=1}^{n}a_{k}\Phi\left(\frac{k}{n}\right)=n^{-\beta}.
\]

This implies

\[
\sum_{k=1}^{n}ka_{k}\left\lfloor \frac{n}{k}\right\rfloor =\sum_{k=1}^{n}\sum_{d\mid k}da_{d}=n^{1-\beta}.
\]

Thus 

\[
\sum_{d\mid n}da_{d}=n^{1-\beta}-\left(n-1\right)^{1-\beta}.
\]

By Möbius inversion, this leads to 

\[
na_{n}=\sum_{d\mid n}\mu\left(\frac{n}{d}\right)\left(d^{1-\beta}-(d-1)^{1-\beta}\right).
\]

Next we have

\[
d^{1-\beta}-(d-1)^{1-\beta}=\left(1-\beta\right)d^{-\beta}+O\left(d^{-1-\beta}\right).
\]

Hence we get
\[
na_{n}=\left(1-\beta\right)\sum_{d\mid n}\mu\left(\frac{n}{d}\right)d^{-\beta}+\sum_{d\mid n}\mu\left(\frac{n}{d}\right)O\left(d^{-1-\beta}\right).
\]

Since $\beta>0$, by previous lemma 2.3 we have

\[
\sum_{d\mid n}\mu\left(\frac{n}{d}\right)d^{-\beta}=O\left(1\right),
\]
and for the second sum we have 
\[
\sum_{d\mid n}\mu\left(\frac{n}{d}\right)O\left(d^{-1-\beta}\right)\ll\sum_{k\leq n}k^{-1-\beta}=O(1).
\]
Thus, we conclude $na_{n}=O(1)$ and it is optimal for $\beta>0$,
since for $p$ prime we have 
\[
pa_{p}=-1+\left(p^{1-\beta}-(p-1)^{1-\beta}\right)\Rightarrow\lim_{p\rightarrow\infty}pa_{p}=-1.
\]
$\square$

\subsubsection{Remark on the HLR criterion}

This criterion is not used in our proof strategy of the Riemann hypothesis
equivalence (Theorem 4.1), based on decomposing $a_{n}$ in the sum
$\sum_{n\leq x}a_{n}$. However, it becomes necessary when approaching
this main theorem through an alternative purely analytic method without
such decomposition. This alternative proof is more direct but relies
more heavily on complex analysis and Mellin convolutions, and will
be developed in a forthcoming paper. \\

The significance of this criterion stems from a fundamental relationship
when $G(n,k)=g(k/n)$ and $g$ is a FGV, as equation (1) creates an
interdependence between the Mellin transform of $g$ and the Dirichlet
series $\sum_{n\geq1}\frac{na_{n}}{n^{s}}$. The Ingham function satisfies
the HLR criterion precisely because of the fundamental theorem of
arithmetic, and this connection manifests through the Euler product
representation of its Mellin transform. This relationship becomes
even more transparent in the context of $L$-functions. Indeed, consider
the generalized Ingham function

\[
\Phi_{u}(x):=\sum_{1\leq k\leq x^{-1}}\frac{u(k)}{k}\Phi\left(kx\right),
\]
where $u=\chi$ is a Dirichlet character. In \cite{key-8} we proved
that its Mellin transform is given by

\[
\Phi_{\chi}^{\star}(z)=\frac{z}{z-1}\zeta\left(1-z\right)L\left(\chi,1-z\right),
\]
and it satisfies the HLR criterion:

\[
\forall\beta\geq0,\ \sum_{k=1}^{n}a_{k}\Phi_{\chi}\left(\frac{k}{n}\right)=n^{-\beta}\Rightarrow\forall\varepsilon>0,\ na_{n}=O\left(n^{\varepsilon}\right).
\]
This criterion holds precisely because $\chi$ is multiplicative,
as reflected in the Euler product representation of $\Phi_{\chi}^{\star}$.
Indeed, if we choose a non-multiplicative function $u$ in the definition
of $\Phi_{u}$, the HLR criterion no longer holds. This structural
connection between the HLR criterion, multiplicative functions and
Euler products could provide new insights into why the Riemann hypothesis
tends to hold for $L$-functions possessing both a functional equation
and an Euler product. This reinforces the view that the Riemann hypothesis
is fundamentally an arithmetic phenomenon.

\section{Four preliminary lemmas}

We establish four lemmas related to equation (1) when $G(n,k)=\Phi\left(\frac{k}{n}\right)$.
These lemmas provide the tools required to establish the asymptotics
central to our Tauberian equivalence.
\begin{itemize}
\item Lemma 3.1 gives estimates of the partial sum of $na_{n}$ when $\beta<0$. 
\item Lemma 3.2 provides expressions for $\sum_{n\geq1}\frac{na_{n}}{n^{s}}$. 
\item Lemma 3.3 gives estimates for partial sums of $\sum_{d\mid n}d^{-\beta}\mu\left(\frac{n}{d}\right)$. 
\item Lemma 3.4 gives estimates for partial sums of $u\star\mu$ when $u(n)=O\left(n^{-\beta-1}\right)$. 
\end{itemize}
We frequently employ the following classical identity

\[
\sum_{n\leq x}\left(u\star v\right)(n)=\sum_{k\leq x}u(k)\sum_{j\leq x/k}v(j).
\]
A proof of which can be found in \cite{key-5} (proposition $4.3$,
p. $177$). 

\subsection{The case $\beta<0$ }

Here, we demonstrate that $\Phi$ is almost a FGV, as it fulfills
the definition of a FGV in the case $\beta<0$. In fact, we will show
something even stronger by introducing an error term in equation (1).

\subsubsection*{Lemma 3.1}

If $\beta<0$ and $0\leq\delta<-\beta$ we have

\[
\sum_{k=1}^{n}a_{k}\Phi\left(\frac{k}{n}\right)=n^{-\beta}+O\left(n^{\delta}\right)\Rightarrow\sum_{n\leq x}a_{n}\sim\frac{1}{\Phi^{\star}\left(\beta\right)}x^{-\beta}\ \left(x\rightarrow\infty\right).
\]

\subsubsection*{Proof of lemma 3.1}

Starting from the classical Möbius inversion formula

\[
\sum_{k=1}^{n}b_{k}\left\lfloor \frac{n}{k}\right\rfloor =T(n)\Rightarrow\sum_{k=1}^{n}b_{k}=\sum_{k=1}^{n}\mu(k)T\left(\left\lfloor \frac{n}{k}\right\rfloor \right),
\]

we obtain

\[
\sum_{k=1}^{n}ka_{k}\left\lfloor \frac{n}{k}\right\rfloor =n^{1-\beta}+O\left(n^{1+\delta}\right)\Rightarrow\sum_{k=1}^{n}ka_{k}=\sum_{k=1}^{n}\mu(k)\left(\left\lfloor \frac{n}{k}\right\rfloor ^{1-\beta}+O\left(\left\lfloor \frac{n}{k}\right\rfloor ^{1+\delta}\right)\right).
\]

Since 

\[
\left\lfloor \frac{n}{k}\right\rfloor ^{1-\beta}=\left(\frac{n}{k}-\left\{ \frac{n}{k}\right\} \right)^{1-\beta}=\left(\frac{n}{k}\right)^{1-\beta}\left(1-\frac{k}{n}\left\{ \frac{n}{k}\right\} \right)^{1-\beta}=\left(\frac{n}{k}\right)^{1-\beta}+\left(\frac{n}{k}\right)^{-\beta}r(n,k)
\]

where $r(n,k)$ is bounded, we get

\[
\sum_{k=1}^{n}ka_{k}=n^{1-\beta}\sum_{k=1}^{n}\frac{\mu(k)}{k^{1-\beta}}+\sum_{k=1}^{n}\mu(k)\left(\frac{n}{k}\right)^{-\beta}r(n,k)+\sum_{k=1}^{n}\mu(k)O\left(\left\lfloor \frac{n}{k}\right\rfloor ^{1+\delta}\right).
\]

Next,

\begin{alignat*}{1}
\left|\sum_{k=1}^{n}\mu(k)\left(\frac{n}{k}\right)^{-\beta}r(n,k)\right| & \ll n^{-\beta}\sum_{k=1}^{n}k^{\beta}=\begin{cases}
O\left(n\right) & \text{for}\quad-1<\beta<0,\\
O\left(n\log n\right) & \text{for}\quad\beta=-1,\\
O\left(n^{-\beta}\right) & \text{for}\quad\beta<-1.
\end{cases}
\end{alignat*}

We also have

\[
\frac{1}{\zeta\left(1-\beta\right)}-\sum_{k=1}^{n}\frac{\mu(k)}{k^{1-\beta}}=\sum_{k\geq n+1}\frac{\mu(k)}{k^{1-\beta}}\ll\sum_{k\geq n+1}\frac{1}{k^{1-\beta}}\ll n^{\beta},
\]

leading to

\[
n^{1-\beta}\sum_{k=1}^{n}\frac{\mu(k)}{k^{1-\beta}}=\frac{n^{1-\beta}}{\zeta\left(1-\beta\right)}+O(n).
\]

And

\[
\sum_{k=1}^{n}\mu(k)O\left(\left\lfloor \frac{n}{k}\right\rfloor ^{1+\delta}\right)\ll\sum_{k=1}^{n}\left\lfloor \frac{n}{k}\right\rfloor ^{1+\delta}\leq n^{1+\delta}\sum_{k=1}^{n}\frac{1}{k^{1+\delta}}=O\left(n^{1+\delta}\log n\right)
\]

since $\delta\geq0$. Thus, we conclude

\[
\sum_{k=1}^{n}ka_{k}=\frac{1}{\zeta\left(1-\beta\right)}n^{1-\beta}+\begin{cases}
O\left(n^{1+\delta}\log n\right) & \text{for}\quad-1<\beta<0,\\
O\left(n^{1+\delta}\log n\right) & \text{for}\quad\beta=-1,\\
O\left(n^{-\beta}\right)+O\left(n^{1+\delta}\log n\right) & \text{for}\quad\beta<-1.
\end{cases}
\]

Therefore, for any $\beta<0$ and $0\leq\delta<-\beta$ we have 

\[
\sum_{k=1}^{n}ka_{k}\sim\frac{1}{\zeta\left(1-\beta\right)}n^{1-\beta}\ \left(n\rightarrow\infty\right).
\]
 Next, by Abel summation, setting $A_{1}(x)=\sum_{n\leq x}na_{n}$,
we obtain

\[
\sum_{n\leq x}a_{n}=\frac{A_{1}(x)}{x}+\int_{1}^{x}\frac{A_{1}(t)}{t^{2}}dt
\]

\[
\sim\frac{1}{\zeta\left(1-\beta\right)}x^{-\beta}+\frac{1}{\zeta\left(1-\beta\right)}\int_{1}^{x}t^{-\beta-1}dt=\frac{1}{\zeta\left(1-\beta\right)}\left(1-\frac{1}{\beta}\right)x^{-\beta}\ \left(x\rightarrow\infty\right)
\]

which finally gives, from the formula for $\Phi^{\star}$ in $3.1,$

\[
\sum_{n\leq x}a_{n}\sim\frac{1}{\Phi^{\star}\left(\beta\right)}x^{-\beta}\ \left(x\rightarrow\infty\right).
\]

$\square$

\subsection{The Dirichlet series for $na_{n}$}

In practice we split the series $\sum_{n\geq1}\frac{na_{n}}{n^{s}}$
into a main term and an error term.

\subsubsection*{Lemma 3.2}

For $\beta>0$ and $\beta\neq1$ we have 

\[
\sum_{n\geq1}\frac{na_{n}}{n^{s}}=\left(1-\beta\right)\frac{\zeta\left(s+\beta\right)}{\zeta\left(s\right)}+\frac{1}{\zeta(s)}\sum_{n\geq1}\frac{u(n)}{n^{s}}
\]
where $u$ is a function satisfying $u(n)=O\left(n^{-\beta-1}\right)$. 

\subsubsection*{Proof of Lemma 3.2}

We have seen in section 2 that

\subsubsection*{
\[
\sum_{k=1}^{n}a_{k}\Phi\left(\frac{k}{n}\right)=n^{-\beta}\Rightarrow na_{n}=\left(1-\beta\right)\sum_{d\mid n}\mu\left(\frac{n}{d}\right)d^{-\beta}+\sum_{d\mid n}\mu\left(\frac{n}{d}\right)O\left(d^{-1-\beta}\right).
\]
}

Therefore,

\[
\sum_{n\geq1}\frac{na_{n}}{n^{s}}=\left(1-\beta\right)\frac{\zeta\left(s+\beta\right)}{\zeta\left(s\right)}+\frac{1}{\zeta(s)}\sum_{n\geq1}\frac{u(n)}{n^{s}}
\]
where $u$ is a function satisfying $u(n)=O\left(n^{-\beta-1}\right)$.
Subsequently, we will use the formula for the partial sums of $na_{n}$

\[
\sum_{n\leq x}na_{n}=\left(1-\beta\right)\sum_{n\leq x}J_{-\beta}(n)+\sum_{n\leq x}(u\star\mu)(n),
\]
where $J_{\lambda}(n)=\sum_{d\mid n}d^{\lambda}\mu\left(\frac{n}{d}\right)$
are the generalized Jordan functions. 

\subsection{Partial sums of generalized Jordan functions}

Here we provide asymptotic formulas for the sum $\sum_{n\leqslant x}J_{-\beta}(n)$.

\subsubsection*{Lemma 3.3}

Suppose $\beta\notin\left\{ 0;1\right\} $, then we have, under the
Riemann hypothesis, when $x\rightarrow\infty$

\[
\left(1-\beta\right)\sum_{n\leqslant x}J_{-\beta}(n)=\begin{cases}
\begin{array}{cc}
\frac{x^{1-\beta}}{\zeta(1-\beta)}+O\left(x^{1/2+\varepsilon}\right), & \ \ 0<\beta<1/2\\
\\
O\left(x^{1/2+\varepsilon}\right), & \beta\geq1/2
\end{array}\end{cases}
\]

\subsubsection*{Proof of lemma 3.3}

\subsubsection*{For $0<\beta<1/2$ }

We have $\sum_{n\geq1}\frac{J_{-\beta}(n)}{n^{s}}=\frac{\zeta\left(s+\beta\right)}{\zeta\left(s\right)}$
which under the Riemann hypothesis
\begin{itemize}
\item Defines a meromorphic function in the half-plane $\Re z>1/2$,
\item Has a single pole at $z=1-\beta>1/2$ with residue $\frac{1}{\zeta\left(1-\beta\right)}$.
\item Has coefficients satisfying $J_{-\beta}(n)=O\left(1\right)$ (cf.
Lemma 2.3).
\end{itemize}
To proceed, let $\delta\in(0,1/2-\beta)$ so that we can apply Perron's
effective formula to the Dirichlet series $\sum_{n\geq1}\frac{J_{-\beta}(n)}{n^{s}}$
( \cite{key-31} II.2, p.133).\\

We integrate over the rectangle with $T=x^{\beta+\delta}$ with the
line of integration $\Re s=1+\frac{1}{\log x}$ and we apply Cauchy's
Residue Theorem to the rectangle 
\[
\left[1-\beta-\delta-iT,1-\beta-\delta+iT,1+\frac{1}{\log x}+iT,1+\frac{1}{\log x}-iT\right].
\]
 The integrand $\frac{\zeta(s+\beta)}{\zeta(s)}\frac{x^{s}}{s}$ has
only one singularity in the rectangle at $s=1-\beta$, which contributes
the residue$\frac{x^{1-\beta}}{\left(1-\beta\right)\zeta(1-\beta)}$.
It remains to bound the integral of $\frac{\zeta(s+\beta)}{\zeta(s)}\frac{x^{s}}{s}$
over three edges of the rectangle: 
\begin{itemize}
\item $s=1-\beta-\delta+it$ with $t\in[-T,T]$,
\item $s=\sigma\pm iT$ with $\sigma\in\left[1-\beta-\delta,1+\frac{1}{\log x}\right]$.
\end{itemize}
To bound $\zeta(s+\beta)$, since the Riemann hypothesis is assumed
to be true, we use the Lindelöf hypothesis. To bound the denominator
$\frac{1}{\zeta(s)}$ we use Theorem $13.23$ in (\cite{key-22},
p.442). This yields, for $|t|\ge1$, $\sigma\ge\frac{1}{2}+\varepsilon$
and some $C>0$

\[
\zeta(\sigma+it)\ll(|t|+C)^{\varepsilon}
\]

\[
\frac{1}{\zeta(\sigma+it)}\ll(|t|+C)^{\varepsilon}
\]

After calculations we obtain
\[
\sum_{n\leqslant x}J_{-\beta}(n)=\frac{x^{1-\beta}}{(1-\beta)\zeta(1-\beta)}+O\left(x^{1-\beta-\delta+\varepsilon}\right),
\]

and taking $\delta=1/2-\beta-\varepsilon$ we get 

\[
\sum_{n\leqslant x}J_{-\beta}(n)=\frac{x^{1-\beta}}{(1-\beta)\zeta(1-\beta)}+O\left(x^{1/2+\varepsilon}\right).
\]

\subsubsection*{For $\beta\geqslant\frac{1}{2}$}

Since the Riemann hypothesis is assumed to be true, we have the Littlewood's
theorem $\sum_{n\leq x}\mu(n)=O\left(x^{1/2+\varepsilon}\right)$
(theorem $3.43$, p.$142$ in \cite{key-5}), therefore we get

\[
\sum_{n\leqslant x}J_{-\beta}(n)=\sum_{k\leqslant x}k^{-\beta}\sum_{j\leqslant x/k}\mu(j)\ll\sum_{k\leqslant x}k^{-\beta}\left(\frac{x}{k}\right)^{1/2+\varepsilon}
\]

\[
\ll x^{1/2+\varepsilon}\sum_{k\leqslant x}k^{-\beta-1/2-\varepsilon}\ll x^{1/2+\varepsilon}.
\]
$\square$

\subsection{Partial sums of a Dirichlet convolution}

Knowledge of the asymptotic behavior of $\sum_{n\leq x}(u\star\mu)(n)$
is necessary when $u(n)=O\left(n^{-\beta-1}\right)$. 

\subsubsection*{Lemma 3.4}

Let $u$ be a function satisfying $u(n)=O\left(n^{-\beta-1}\right),$
then, if $\beta\geq-1/2$, we have

\[
\sum_{n\leq x}(u\star\mu)(n)\ll x^{1/2+\varepsilon}
\]

\paragraph*{Proof of lemma 3.4}

Since we assume that the Riemann hypothesis is true, we have

\[
\sum_{n\leq x}(u\star\mu)(n)=\sum_{k\leqslant x}u(k)\sum_{j\leqslant x/k}\mu(j)\ll\sum_{k\leqslant x}\left|u(k)\right|\left(\frac{x}{k}\right)^{1/2+\varepsilon},
\]
and therefore,

\[
\sum_{n\leq x}(u\star\mu)(n)\ll x^{1/2+\varepsilon}\sum_{k\leqslant x}k^{-\beta-3/2}.
\]
Since $\beta\geq-1/2$ we have $\beta+3/2\geq1$, thus we get

\[
\sum_{n\leq x}(u\star\mu)(n)\ll x^{1/2+\varepsilon}\log x\ll x^{1/2+\varepsilon}.
\]
$\square$

\section{Tauberian equivalence of the Riemann hypothesis}

The central result of this paper is a Tauberian equivalence of the
Riemann hypothesis, which is the subject of the following theorem. 

\subsection{Equivalence theorem }

The following two propositions are equivalent:
\begin{description}
\item [{(i)}] The Riemann hypothesis is true.
\item [{(ii)}] $\Phi$ is a function of good variation with index of regularity
$\alpha\left(\Phi\right)=\frac{1}{2}$.
\end{description}

\subsection{Proof of the Tauberian equivalence theorem 4.1}

We first demonstrate the simpler implication that $(ii)$ implies
$(i)$. 

\subsubsection{Proof of (ii) implies (i)}

Assume that $\Phi$ is a function of good variation (FGV) with index
$\alpha\left(\Phi\right)=\frac{1}{2}$. Setting $\beta=1>\frac{1}{2}$,
we have, by definition of a FGV

\[
\sum_{k=1}^{n}a_{k}\Phi\left(\frac{k}{n}\right)=n^{-1}\Rightarrow\sum_{k=1}^{n}a_{k}=O\left(n^{-1/2+\varepsilon}\right).
\]
Since for all $n\geq1,$
\[
\sum_{k=1}^{n}a_{k}\Phi\left(\frac{k}{n}\right)=n^{-1}\Longleftrightarrow\sum_{k=1}^{n}ka_{k}\left\lfloor \frac{n}{k}\right\rfloor =1\Longleftrightarrow na_{n}=\mu_{n},
\]
where $\mu$ is the Möbius function (this is due to the Meissel identity
$\sum_{k=1}^{n}\mu_{k}\left\lfloor \frac{n}{k}\right\rfloor =1$ \cite{key-2}),
it follows that

\[
\sum_{k=1}^{n}\frac{\mu_{k}}{k}=O\left(n^{-1/2+\varepsilon}\right)\Longleftrightarrow\sum_{k=1}^{n}\mu_{k}=O\left(n^{1/2+\varepsilon}\right),
\]
using Abel summation, which is equivalent to the Riemann hypothesis.$\square$

\subsubsection{Proof of (i) implies (ii) }

To address this implication, we analyze the cases for $\beta$ in
the equation

\[
\sum_{k=1}^{n}a_{k}\Phi\left(\frac{k}{n}\right)=n^{-\beta}.
\]
We first examine the cases $\beta\in\left\{ 0;1\right\} $. 

\paragraph{Case $\beta=0$}

We have $a_{1}=1$ and $a_{n}=0$ for $n\geq2$, hence we obtain

\[
\forall n\geq1,\ \sum_{k=1}^{n}a_{k}\Phi\left(\frac{k}{n}\right)=\sum_{k=1}^{n}a_{k}=1.
\]
Then,

\[
\Phi^{\star}(0)=\lim_{z\rightarrow0}\frac{z}{z-1}\zeta\left(1-z\right)=1.
\]
Thus we can write

\begin{equation}
\beta=0\Rightarrow\sum_{n\leq x}a_{n}\sim\frac{1}{\Phi^{\star}(\beta)}x^{-\beta}\ \left(x\rightarrow\infty\right).
\end{equation}

\paragraph{Case $\beta=1$ }

From $4.2.1$, we know that $a_{n}=\frac{\mu(n)}{n}$ . Assuming the
Riemann hypothesis, we obtain

\[
\sum_{n\leq x}\frac{\mu(n)}{n}=O\left(x^{-1/2+\varepsilon}\right).
\]
Thus, we can write 

\begin{equation}
\beta=1\Rightarrow\sum_{n\leq x}a_{n}=O\left(x^{-1/2+\varepsilon}\right).
\end{equation}

\paragraph{Case $\beta<0$ }

By lemma $3.1$, we have the following relationship:

\begin{equation}
\beta<0\Rightarrow\sum_{n\leq x}a_{n}\sim\frac{1}{\Phi^{\star}\left(\beta\right)}x^{-\beta}\ \left(x\rightarrow\infty\right).
\end{equation}

\subsubsection*{Case $0<\beta<1/2$ }

According to lemma $3.3$, 

\[
\left(1-\beta\right)\sum_{n\leq x}J_{-\beta}(n)=\frac{x^{1-\beta}}{\zeta(1-\beta)}+O\left(x^{1/2+\varepsilon}\right),
\]
and according to lemma $3.4$, 

\[
\sum_{n\leq x}u\star\mu(n)=O\left(x^{1/2+\varepsilon}\right).
\]
Thus, we obtain

\[
\sum_{n\leq x}na_{n}=\frac{x^{1-\beta}}{\zeta(1-\beta)}+O\left(x^{1/2+\varepsilon}\right).
\]
Applying Abel summation, we conclude

\begin{equation}
0<\beta<1/2\Rightarrow\sum_{n\leq x}a_{n}\sim\frac{1}{\Phi^{\star}\left(\beta\right)}x^{-\beta}\ \left(x\rightarrow\infty\right).
\end{equation}

\subsubsection*{Case $\beta\protect\geq1/2$ }

According to lemma $3.3$, 

\[
\left(1-\beta\right)\sum_{n\leq x}J_{-\beta}(n)=O\left(x^{1/2+\varepsilon}\right).
\]
And according to lemma $3.4$,

\[
\sum_{n\leq x}u\star\mu(n)=O\left(x^{1/2+\varepsilon}\right).
\]
Thus, we obtain

\[
\sum_{n\leq x}na_{n}=O\left(x^{1/2+\varepsilon}\right).
\]
Applying Abel summation, we conclude

\begin{equation}
\beta\geq1/2\Rightarrow\sum_{n\leq x}a_{n}=O\left(x^{-1/2+\varepsilon}\right)\ \left(x\rightarrow\infty\right).
\end{equation}

\subsubsection*{Summary of the proof of (i) implies (ii) }

Let's summarize: formulas $(8),(9),(10),(11)$ and $\left(12\right)$
lead to the following conclusions:

\[
\beta<\frac{1}{2}\Rightarrow\sum_{n\leq x}a_{n}\sim\frac{1}{\Phi^{\star}\left(\beta\right)}x^{-\beta}\ \left(x\rightarrow\infty\right),
\]

\[
\beta\geq1/2\Rightarrow\sum_{n\leq x}a_{n}\ll x^{-1/2+\varepsilon}\ \left(x\rightarrow\infty\right).
\]
These asymptotic behaviors precisely match the definition of a function
of good variation (FGV) given in the introduction, and thus, the Ingham
function $\Phi$ is a FGV with regularity index $\alpha\left(\Phi\right)=\frac{1}{2}$.$\square$

\subsubsection{Verification for $\beta=\infty$}

We verify that our results hold when $\beta=\infty>1/2$. In this
case, with $a_{1}=1$, for $n\geq2$ equation $(1)$ simplifies to

\[
\sum_{k=1}^{n}a_{k}\Phi\left(\frac{k}{n}\right)=0
\]
 Using the formula
\[
na_{n}=\mu(n)+\sum_{d\mid n,d\geq2}\mu\left(\frac{n}{d}\right)\left(d^{1-\beta}-(d-1)^{1-\beta}\right),
\]
and noting that $d^{1-\beta}=0$ for $d\geq2$ and $(d-1)^{1-\beta}=0$
for $d\geq3$, we find, for $k\geq1$ 

\[
a_{2k}=\frac{\mu\left(2k\right)-\mu\left(k\right)}{2k}\quad\text{and}\quad a_{2k+1}=\frac{\mu\left(2k+1\right)}{2k+1}.
\]
It can be noted that the sequence $na_{n}$ is referenced as sequence
\href{https://oeis.org/A092673}{A092673} in \cite{key-24}. Define
$M(x):=\sum_{n\leq x}\mu(n)$ and $A_{1}(x):=\sum_{n\leq x}na_{n}$.
Then it is easy to see

\[
M(x)=A_{1}(x)+A_{1}\left(\frac{x}{2}\right)+A_{1}\left(\frac{x}{2^{2}}\right)+...
\]
with $A_{1}(0)=0$. Hence we get

\[
A_{1}(x)=M(x)-M\left(\frac{x}{2}\right).
\]
Assuming the Riemann hypothesis is true, we have $M(x)=O\left(x^{1/2+\varepsilon}\right)$,
which implies
\[
A_{1}(x)=O\left(x^{1/2+\varepsilon}\right)\Rightarrow\sum_{n\leq x}a_{n}=O\left(x^{-1/2+\varepsilon}\right).
\]
If, on the other hand, $\Phi$ is a FGV with regularity index $\alpha\left(\Phi\right)=\frac{1}{2}$,
then

\[
\sum_{k=1}^{n}a_{k}\Phi\left(\frac{k}{n}\right)=0\Rightarrow\sum_{n\leq x}a_{n}=O\left(x^{-1/2+\varepsilon}\right).
\]
Therefore we obtain

\[
A_{1}(x)=O\left(x^{1/2+\varepsilon}\right)
\]
and consequently

\[
\left|M(x)\right|\leq\sum_{n\geq0}\left|A_{1}\left(\frac{x}{2^{n}}\right)\right|\ll\sum_{n\geq0}\left(\frac{x}{2^{n}}\right)^{1/2+\varepsilon}\ll x^{1/2+\varepsilon},
\]
confirming that the Riemann hypothesis is true.\\

This limiting case not only confirms the consistency of our Tauberian
equivalence of the Riemann hypothesis but also highlights the Möbius
function's pivotal role, showcasing the deep interplay between arithmetic
and analysis in our framework.\\

We now establish a corollary of Theorem 4.1 that connects the regularity
index of the Ingham function to the non-trivial zeros of the zeta
function.

\subsection{Corollary of the equivalence theorem}

This corollary extends our main result by providing a precise correspondence
between the zero-free regions of the zeta function and the regularity
index of the Ingham function.

\subsubsection*{Corollary 4.3}

Let $1/2\leq\sigma\leq1$. Then the following propositions are equivalent:
\begin{description}
\item [{(i)}] The Riemann zeta function has no zeros in the half-plane
$\Re z>\sigma$,
\item [{(ii)}] $\Phi$ is a FGV with regularity index $1-\sigma\leq\alpha\left(\Phi\right)\leq1/2$.
\end{description}

\subsubsection*{Proof of corollary $4.3$}

The proof follows the same steps as that of Theorem $4.1$ by adapting
Lemmas $3.3$ and $3.4$.$\square$

\subsection{Recovering the prime number theorem}

An interesting observation is that the mere knowledge that $\Phi$
is a FGV with strictly positive regularity index (assuming we only
know that $\alpha\left(\Phi\right)>0$, without using Corollary 4.3),
and with no information about the zeros of the zeta function, is sufficient
to recover the prime number theorem.

\subsubsection*{Proof}

Assuming $\Phi$ is a FGV with regularity $\alpha\left(\Phi\right)>0$,
since $\infty>\alpha\left(\Phi\right)$, we get from FGV properties
for $n\geq2$

\[
\sum_{k=1}^{n}a_{k}\Phi\left(\frac{k}{n}\right)=0\Rightarrow\sum_{n\leq x}a_{n}=O\left(x^{-\alpha\left(\Phi\right)+\varepsilon}\right)=o\left(1\right)
\]
and thus by Abel summation

\[
\sum_{n\leq x}na_{n}=o\left(x\right).
\]
Next the calculations in $4.2.3$ allow us to say that

\[
\sum_{n\leq x}na_{n}=M(x)-M\left(\frac{x}{2}\right)=o\left(x\right)\Rightarrow M(x)=\sum_{k\geq0}M\left(\frac{x}{2^{k}}\right)-M\left(\frac{x}{2^{k+1}}\right)
\]

\[
\ll\sum_{k\geq0}o(x/2^{k})=o(x),
\]
which is equivalent to the prime number theorem (see for instance
\cite{key-2}).\\

\subsubsection*{Remark}

Lemma 3.1 already shows that $\alpha\left(\Phi\right)\geq0$. However,
proving $\alpha\left(\Phi\right)>0$ would imply, by Corollary 4.3,
that $\zeta\left(s\right)$ has no zeros for $1-\alpha\left(\Phi\right)<\Re s<1$,
which is still an open problem in analytic number theory. Therefore,
one should seek non-analytic methods to achieve this goal and there
might exist a general rule to determine whether a given function $g$
is an FGV with strictly positive regularity index. This is an interesting
research direction in its own right within RAF theory.\\

\section{The birth of combinatorial Tauberian theory}

Although the primary focus of this paper is the Tauberian equivalence
of the Riemann hypothesis, this section demonstrates how RAF theory
extends naturally to combinatorial number-theoretic problems. This
extension leverages the intrinsic structure of RAF theory, which unifies
discrete and continuous aspects of number theory. By bridging analytic
and combinatorial perspectives, it introduces the foundations of what
we term combinatorial Tauberian theory. This common branch between
Tauberian and combinatorial theories (specially combinatorial number
theory) appears to be new, as terms such as ``combinatorial'', ``enumeration'',
and ``counting'' yield no significant results in \cite{key-20}.

\subsection*{Definitions}

\subsubsection*{Combinatorial kernel}

A function $C(n,k)$ is called a combinatorial kernel of natural numbers
if there exist at least two sequences of rational numbers $u_{n}$
and $v_{n}$ such that:

\[
\sum_{k=1}^{n}u_{k}C(n,k)=\sum_{k=1}^{n}v_{k}
\]
and where $V(n):=\sum_{k=1}^{n}v_{k}$ is an integer value counting
a specific quantity.\\

Two types of problems arise from this definition, reflecting fundamental
questions in enumerative combinatorics and asymptotic number theory.

\subsubsection*{Combinatorial Abelian problems}

Given a combinatorial kernel $C(n,k)$ and known $u_{n}$, find an
exact formula or estimate an asymptotic formula for $\sum_{k=1}^{n}u_{k}C(n,k)$.

\subsubsection*{Combinatorial Tauberian problems}

Given a combinatorial kernel $C(n,k)$, known $V(n)$, find an exact
formula or estimate an asymptotic formula for $U(n)=\sum_{k=1}^{n}u_{k}$.\\

Our primary interest lies in combinatorial kernels that facilitate
counting within number theory rather than in combinatorics in general.
We therefore exclude kernels such as binomial coefficients $C(n,k)={n \choose k}$
or Stirling numbers, as these offer limited Tauberian insights and
can be handled through existing inversion formulas (inverse binomial
or Stirling transform) without requiring RAF theory. While combinatorial
kernels are not generally regular arithmetic functions themselves,
they can be transformed into RAF by appropriate scaling factors, allowing
the application of our theory. 

From this point forward, we will focus on the specific kernel $C(n,k)=\left\lfloor \frac{n}{k}\right\rfloor $
which, as we shall see, yields number-theoretic counting problems
with Tauberian implications within RAF theory.

\paragraph*{Example of number theoretic combinatorial Abelian problem}

$C(n,k)=\left\lfloor \frac{n}{k}\right\rfloor $ is a combinatorial
kernel because if we take $u_{n}=1$, we have 

\[
\sum_{k=1}^{n}u_{k}C(n,k)=\sum_{k=1}^{n}\left\lfloor \frac{n}{k}\right\rfloor =\sum_{k=1}^{n}\tau(k)
\]
where $\tau$ is the divisor-counting function. Therefore, the sum
$\sum_{k=1}^{n}u_{k}C(n,k)$ counts the total number of divisors of
all integers up to $n$. Finding a precise asymptotic formula for
$\sum_{k=1}^{n}\tau(k)$ is the famous Dirichlet divisor problem,
and we are facing here a highly non-obvious combinatorial Abelian
problem (\cite{key-5}, p.238, p.307).

\paragraph*{Remark on Tauberian conditions}

Unlike in classical Tauberian theory, we do not impose conditions
on $u_{n}$ (such as growth conditions) because in combinatorial Tauberian
theory, the sequence $u_{n}$ is uniquely determined through the recurrence
relation $\sum_{k=1}^{n}u_{k}C(n,k)=\sum_{k=1}^{n}v_{k}$ which forces
the behavior of the sequence $u_{n}$. This is somewhat similar to
the HLR criterion discussed in Section 2.

\subsection{Generalization of Equation (1)}

Equation (1) can be generalized in several ways to address additional
problems in number theory. As seen in Lemma 3.1, where we could add
an error term to the RHS, we can also consider regularly varying functions
with index $-\beta$ (in the sense of Karamata, see \cite{key-20},
chapter IV, pp. 177-232), namely:

\[
\sum_{k=1}^{n}a_{k}G(n,k)=n^{-\beta}L_{0}(n),
\]
where $L_{0}$ is slowly varying (i.e $\forall\lambda>0,\ \lim_{x\rightarrow\infty}\frac{L_{0}\left(\lambda x\right)}{L_{0}(x)}=1)$,
increasing on $[1,\infty)$ and satisfying $0\leq L_{0}(n+1)-L_{0}(n)<C$.
We can show this general result (details will appear in a forthcoming
paper). If $G$ is a RAF with regularity index $\alpha\left(G\right)$,
then we have
\[
\beta<\alpha\left(G\right)\ \Rightarrow\ \sum_{n\leq x}a_{n}\sim\frac{1}{G^{\star}\left(\beta\right)}x^{-\beta}L_{0}(x)\ \left(x\rightarrow\infty\right)
\]

\[
\beta\geq\alpha\left(G\right)\ \Rightarrow\sum_{n\leq x}a_{n}\ll x^{-1/2+\varepsilon}\ \left(x\rightarrow\infty\right).
\]

Here, we focus on the case of Ingham's function and we prove the following
partial theorem which is sufficient for our purposes.

\subsubsection*{Theorem 5.1}

Suppose $a_{1}=1$, and for $n\geq2$: 

\[
\sum_{k=1}^{n}a_{k}\Phi\left(\frac{k}{n}\right)=n^{-\beta}L_{0}(n),
\]
then, under the Riemann hypothesis, we have:

\[
\beta\geq1/2\ \Rightarrow\sum_{n\leq x}a_{n}\ll x^{-1/2+\varepsilon}\ \left(x\rightarrow\infty\right).
\]

\subsubsection*{Proof of theorem 5.1}

We have

\[
\sum_{k=1}^{n}ka_{k}\left\lfloor \frac{n}{k}\right\rfloor =\sum_{k=1}^{n}\sum_{d\mid k}da_{d}=n^{1-\beta}L_{0}(n).
\]
hence we get

\[
\sum_{d\mid n}da_{d}=n^{1-\beta}L_{0}(n)-\left(n-1\right)^{1-\beta}L_{0}(n-1).
\]
By Möbius inversion, this leads to 

\[
na_{n}=\sum_{d\mid n}\mu\left(\frac{n}{d}\right)\left(d^{1-\beta}L_{0}(d)-(d-1)^{1-\beta}L_{0}(d-1)\right).
\]
Next, since we have $(d-1)^{1-\beta}=d^{1-\beta}-\left(1-\beta\right)d^{-\beta}+O\left(d^{-1-\beta}\right)$,
letting
\begin{itemize}
\item $u(n)=\left(1-\beta\right)n^{-\beta}L_{0}(n-1)$,
\item $v(n)=O\left(n^{-1-\beta}L_{0}(n)\right),$ 
\item $w(n)=n^{1-\beta}\left(L_{0}(n)-L_{0}(n-1)\right),$ 
\end{itemize}
we get

\[
na_{n}=\left(u\star\mu\right)(n)+\left(v\star\mu\right)(n)+\left(w\star\mu\right)(n).
\]
Under the Riemann hypothesis we then have:

\[
\sum_{n\leq x}u\star\mu(n)\ll\sum_{k\leq x}u(k)(x/k)^{1/2+\varepsilon}
\]

\[
\ll x^{1/2+\varepsilon}\sum_{k\leq x}k^{-\beta}L_{0}(k-1)k^{-1/2-\varepsilon}\ll x^{1/2+\varepsilon}\sum_{k\leq x}k^{-\beta-1/2+\varepsilon}\ll x^{1/2+\varepsilon}
\]
 since $L_{0}(n)=O\left(n^{\varepsilon}\right)$ and $\beta\geq1/2$.
Similarly we have 

\[
\sum_{n\leq x}v\star\mu(n)\ll\sum_{k\leq x}v(k)(x/k)^{1/2+\varepsilon}
\]

\[
\ll x^{1/2+\varepsilon}\sum_{k\leq x}k^{-1-\beta}L_{0}(k)k^{-1/2-\varepsilon}\ll x^{1/2+\varepsilon}\sum_{k\leq x}k^{-3/2+\varepsilon}\ll x^{1/2+\varepsilon}
\]
Next

\[
\sum_{n\leq x}w\star\mu(n)=\sum_{k\leq x}w(k)\sum_{j\leq x/k}\mu(j)
\]

\[
\ll\sum_{k\leq x}k^{1-\beta}\left(L_{0}(k)-L_{0}(k-1)\right)(x/k)^{1/2+\varepsilon}
\]

\[
\ll x^{1/2+\varepsilon}\sum_{k\leq x}k^{1/2-\beta-\varepsilon}\left(L_{0}(k)-L_{0}(k-1)\right).
\]
Now, using Abel summation, we get

\[
\sum_{k\leq x}k^{1/2-\beta-\varepsilon}\left(L_{0}(k)-L_{0}(k-1)\right)=x^{1/2-\beta-\varepsilon}L_{0}(x)-\sum_{k\leq x}\left(k^{1/2-\beta-\varepsilon}-(k-1)^{1/2-\beta-\varepsilon}\right)L_{0}(k)\ll x^{1/2-\beta+\varepsilon},
\]
since $L_{0}(k)=O\left(k^{\varepsilon}\right)$ and $\left(k^{1/2-\beta-\varepsilon}-(k-1)^{1/2-\beta-\varepsilon}\right)L_{0}(k)\ll k^{-1/2-\beta+\varepsilon}$.
Thus we get 

\[
\sum_{n\leq x}w\star\mu(n)\ll x^{1-\beta+\varepsilon},
\]
yielding, since $\beta\geq\frac{1}{2}$,

\[
\sum_{n\leq x}na_{n}=\sum_{n\leq x}u\star\mu(n)+\sum_{n\leq x}v\star\mu(n)+\sum_{n\leq x}w\star\mu(n)\ll n^{1/2+\varepsilon}.
\]
$\square$

We can now illustrate how the combinatorial properties of the Möbius
and floor functions, central to many counting problems in number theory,
provide a natural playground allowing us to use RAF theory in relation
with combinatorial number theory. 

\subsection{Combinatorial properties of the Floor and Möbius functions}

Together, the floor function and Möbius function generate many combinatorial
formulas counting integers with specific arithmetic properties. Here
are some examples:

\subsubsection*{Counting coprime $m$-tuples}

$\sum_{k=1}^{n}\mu(k)\left\lfloor \frac{n}{k}\right\rfloor ^{2}$
counts the number of pairs 

\[
\left\{ \left(x_{1},x_{2}\right),\ 1\leq x_{1},x_{2}\leq n,\ \gcd\left(x_{1},x_{2}\right)=1\right\} 
\]
which is sequence \href{https://oeis.org/A018805}{A018805} in \cite{key-24}.
More generally, for $m\geq2$, the number $\sum_{k=1}^{n}\mu(k)\left\lfloor \frac{n}{k}\right\rfloor ^{m}$
counts unordered $m$-tuples:

\[
\left\{ \left(x_{i}\right)_{i=1,..,m},\ 1\leq x_{1},...,x_{m}\leq n,\ \gcd\left(x_{1},...,x_{m}\right)=1\right\} 
\]

There are several other sequences in \cite{key-24} : \href{https://oeis.org/A071778}{A071778},\href{https://oeis.org/A082540}{A082540},\href{https://oeis.org/A082544}{A082544}
for $m=3,4,5,...$ respectively.\\

Note that Meissel's identity $\sum_{k=1}^{n}\mu(k)\left\lfloor \frac{n}{k}\right\rfloor =1$
correctly counts numbers such that $\gcd(x_{1},x_{1})=1$, since only
$x_{1}=1$ works, making the formula valid for $m\geq1$.

\subsubsection*{Counting $p$-free integers}

$\sum_{k=1}^{n}\mu(k)\left\lfloor \frac{n}{k^{2}}\right\rfloor $
counts squarefree integers up to $n$ which is sequence \href{https://oeis.org/A013928}{A013928}
in \cite{key-24}. This formula is particularly useful for asymptotic
estimates \cite{key-9}.\\

More generally, for prime $p\geq2$, the number $\sum_{k=1}^{n}\mu(k)\left\lfloor \frac{n}{k^{p}}\right\rfloor $
counts $p$-free integers up to $n$ and the sequence \href{https://oeis.org/A060431}{A060431}
gives the number of cubefree integers up to $n$.

\subsubsection*{Counting integers with determined prime factors}

$-\sum_{k=1}^{n}\mu(2k)\left\lfloor \frac{n}{k}\right\rfloor $ counts
powers of two up to $n$ which is sequence \href{https://oeis.org/A070939}{A070939}
in \cite{key-24}.\\

More generally, for prime $p\geq2$, the number $-\sum_{k=1}^{n}\mu(pk)\left\lfloor \frac{n}{k}\right\rfloor $
counts powers of $p$ up to $n$ and it's straightforward to show
that 
\[
-\sum_{k=1}^{n}\mu(pk)\left\lfloor \frac{n}{k}\right\rfloor =1+\left\lfloor \frac{\log n}{\log p}\right\rfloor .
\]

Further generalizing, if $P=\prod_{i=1}^{m}p_{i}$, where $p_{i}$
are distinct primes, then

\[
(-1)^{m}\sum_{k=1}^{n}\mu(Pk)\left\lfloor \frac{n}{k}\right\rfloor 
\]

counts integers up to $n$ containing only $p_{i}$ in their prime
factorization. There are several sequences in \cite{key-24} :\\

\href{https://oeis.org/A081604}{A081604},\href{https://oeis.org/A110592}{A110592},\href{https://oeis.org/A071521}{A071521},\href{https://oeis.org/A071520}{A071520},\href{https://oeis.org/A071604}{A071604},\href{https://oeis.org/A112751}{A112751}
\\

for $P=3,5,6,30,210,15$ respectively.

\subsection{Connection between RAF theory and combinatorial number theory}

While the first two examples above are of no immediate relevance to
our discussion, the third example, however, highlights how RAF theory
naturally interacts with combinatorial number theory, particularly
through the combinatorial kernel $C(n,k)=\left\lfloor \frac{n}{k}\right\rfloor $.
For instance, let us consider the problem of estimating the asymptotic
behavior of the sum $\sum_{n\leq x}\mu\left(6n\right)$ under the
Riemann hypothesis. Below, we present the classical analytic approach
followed by the Tauberian combinatorial one.

\subsubsection{Pure analytical approach}

It is easy to see that the Dirichlet series for $\mu\left(6n\right)$
is given by

\[
F(s)=\sum_{n\geq1}\frac{\mu\left(6n\right)}{n^{s}}=\frac{1}{\left(1-2^{-s}\right)\left(1-3^{-s}\right)\zeta\left(s\right)}=\prod_{p\geq5}\frac{1}{1-p^{-s}}
\]

The zeros of $\left(1-2^{-s}\right)\left(1-3^{-s}\right)$ lie on
the line $x=0$ and, under the Riemann hypothesis, $F$'s rightmost
singularities are given by zeta's zeros which are on the critical
line. Thus a classical analytical calculation yields (details omitted):

\[
\sum_{n\leq x}\mu\left(6n\right)\ll x^{1/2+\varepsilon}
\]

While this approach highlights the analytic properties of the Möbius
function and the zeta function, RAF theory offers an alternative pathway
rooted in its combinatorial structure.

\subsubsection{Combinatorial Tauberian approach}

RAF theory enables us to frame this result within a more arithmetical
framework where the good variation properties of the Ingham function
provide a direct connection to the Riemann hypothesis. As seen in
5.2, the formula
\[
\sum_{k=1}^{n}\mu\left(6k\right)\left\lfloor \frac{n}{k}\right\rfloor 
\]
counts $3$-smooth integers up to $n$ which is sequence \href{https://oeis.org/A071521}{A071521}
in \cite{key-24}. \\

Ramanujan \cite{key-14} showed this number is a function $L_{0}(n)$
such that:

\[
L_{0}(n)\sim\frac{\log\left(2n\right)\log\left(3n\right)}{2\log(2)\log(3)}\ \left(n\rightarrow\infty\right).
\]
Since $L_{0}$ is a slowly varying function and we trivially have
$0\leq L_{0}(n+1)-L_{0}(n)\leq1$, from the equation:

\[
\sum_{k=1}^{n}\frac{\mu\left(6k\right)}{k}\Phi\left(\frac{k}{n}\right)=n^{-1}L_{0}\left(n\right),
\]
using the theorem 5.1 (here we have $1>1/2$), we get

\[
\sum_{k=1}^{n}\frac{\mu\left(6k\right)}{k}\ll n^{-1/2+\varepsilon},
\]
which by Abel summation recovers:

\[
\sum_{n\leq x}\mu\left(6n\right)\ll x^{1/2+\varepsilon}.
\]

\subsubsection{More exotic enumeration}

We can show in the same combinatorial Tauberian way (cf. 5.3.2) that

\[
\sum_{n\leq x}(-1)^{n-1}\mu\left(n\right)\ll x^{1/2+\varepsilon}.
\]

Indeed, it is easy to see that we have

\[
\sum_{k=1}^{n}(-1)^{k-1}\mu\left(k\right)\left\lfloor \frac{n}{k}\right\rfloor =1+2\left\lfloor \frac{\log n}{\log2}\right\rfloor ,
\]
which is sequence \href{https://oeis.org/A129972}{A129972} in \cite{key-24}
which counts the number of bits needed to write $n$ using Elias gamma
coding \cite{key-26}.\\

It would be interesting to find other examples in combinatorial number
theory where RAF theory could be used to estimate non-trivial asymptotic
formulas for partial sums of significant and diverse arithmetic functions,
using various RAF $G$ not necessarily linked to the Ingham function.\\

While the examples provided here demonstrate the potential of RAF
theory to address combinatorial questions, they also suggest deeper
structural connections. The following section extends this perspective,
unveiling an unexpected link between number theory and algebraic geometry
through RAF theory. This connection not only enriches our understanding
of regular arithmetic functions but also suggests new approaches to
fundamental problems like the Riemann hypothesis and its generalizations.

\section{Perspective}

The proof of our Tauberian equivalence of the Riemann hypothesis is
both natural and straightforward, relying on Möbius inversion applied
to weighted sums---also an indispensable technique in elementary
proofs of the prime number theorem (see for instance \cite{key-2}).
We first employ this technique in Lemma 3.1 to sums involving the
floor function, then apply it again in Lemma 3.2 to the Dirichlet
convolution, which ultimately leads to the key estimates presented
in Lemmas 3.3 and 3.4.. This is not surprising, as Tauberian theorems
are fundamentally inversion formulas for weighted sums.

What makes this equivalence particularly promising is its potential
to be integrated into a broader conceptual framework. By introducing
an initially ``unmotivated'' new parameter\footnote{This is a reference to Freeman Dyson's remark in an interview where,
discussing his paper on ferromagnet, he directly alludes to Littlewood's
introduction of an \textquotedbl unmotivated parameter\textquotedbl{}
that directly inspired him: \textquotedbl So that was the big Tauberian
theorem which Littlewood proved, and this way of proving it was again
a tour de force done by invoking a completely unmotivated new parameter
which turned out to be the key to the proof.\textquotedbl} (here, a function $f$), we lay the groundwork for a more abstract
interpretation of the Riemann hypothesis, viewing the regularity index
as an invariant understood through the notions of stability and balance
of a RAF with respect to $f$. This perspective opens the door to
future developments, which we intend to explore in a forthcoming paper
and which we briefly outline here.

\subsection{Stability and balance of a RAF }

Given a RAF $G$ and a function $f,$ we consider a generalization
of equation (1), namely

\[
\sum_{k=1}^{n}a_{k}G\left(f(n),f(k)\right)=f(n)^{-\beta}.
\]

This equation motivates the introduction of two key notions: stability
and balance of a RAF $G$ with respect to a function $f$. These notions,
in turn, necessitate extending the definitions of the arithmetic Mellin
transform and the regularity index.

\subsubsection*{Arithmetic Mellin transform with respect to $f$}

Let $f$ be a positive increasing function. We say that $G$ has an
arithmetic Mellin transform with respect to $f$ if, for $\Re z<0$,
the following limit exists: 
\[
G_{f}^{\star}(z):=\lim_{n\rightarrow\infty}f(n)^{z}\sum_{k=1}^{n}\left(f(k)^{-z}-f(k-1)^{-z}\right)G\left(f(n),f(k)\right).
\]

\subsubsection*{Regularity index of a RAF with respect to \textmd{$f$}}

For a RAF $G$ and a positive function $f$ increasing to infinity,
we define the regularity Index of $G$ with respect to $f$, denoted
$\alpha_{f}(G)$, as the supremum of all $\beta$ such that

\[
\sum_{k=1}^{n}a_{k}G(f(n),f(k))=f(n)^{-\beta}\Rightarrow\sum_{k=1}^{n}a_{k}\sim\frac{1}{G_{f}^{\star}(\beta)}f(n)^{-\beta}\ \left(n\rightarrow\infty\right)
\]

\subsubsection*{Stability of a RAF with respect to \textmd{$f$}}

A RAF $G$ is said to be stable with respect to $f$ if $G$ possesses
a regularity index with respect to $f$.

\subsubsection*{Balance of a RAF with respect to \textmd{$f$}}

A RAF $G$ is said to be balanced with respect to $f$ if it is stable
with respect to $f$ and, moreover, $\alpha_{f}(G)=\alpha(G)$.

\subsection{Balance of the Ingham function}

The Ingham function exhibits remarkable properties in our framework.
It is balanced for functions such as $f(x)=x^{r}\left(1+o\left(1\right)\right)$
($0<r\leq1$) since here we can show $\Phi_{f}^{\star}(z)=\Phi^{\star}(z)$
and $\alpha_{f}\left(\Phi\right)=\alpha\left(\Phi\right)$ unconditionally
and, on the other hand, it is stable with respect to functions with
exponential growth. Particularly if $f(x)=q^{x}+1$ with $q\geq2$
integer, we obtain the arithmetic Mellin transform of the Ingham function
with respect to $f$:
\[
\Phi_{f}^{\star}(z)=\frac{q^{2z}-2q^{z}+q}{q-q^{z}}
\]

which zeros lie on the critical line $x=1/2$, yielding $\alpha_{f}(\Phi)=1/2$.
Hence proving that the Ingham function is balanced (and not merely
stable) with respect to $f(x)=2^{x}+1$ would be equivalent to proving
the Riemann hypothesis.\\

Moreover, we observe a striking connection between the Ingham function
and the Hasse zeta function of an elliptic curve over a finite field
\cite{key-15}. Specifically, letting $a=2$ and $T=q^{-z}$, we can
express $\Phi_{f}^{\star}(z)$ given above as follows

\[
\Phi_{f}^{\star}(z)=\left(1-\frac{1}{T}\right)Z\left(C/\mathbb{F}_{q},T\right)
\]
where $Z\left(C/\mathbb{F}_{q},T\right)$ is the Hasse zeta function 

\[
Z\left(C/\mathbb{F}_{q},T\right)=\frac{qT^{2}-aT+1}{\left(1-T\right)\left(1-qT\right)}.
\]

This connection intriguingly links the Riemann zeta function and the
zeta function of an elliptic curve over a finite field---two entities
traditionally considered to belong to distinct mathematical domains
(analytic number theory vs algebraic geometry). This surprising bridge
hints at a unified arithmetic framework, reinforcing the potential
of RAF theory to uncover deep structural truths underlying the Riemann
hypothesis.\\

The findings presented here are the starting point of ongoing work
aimed at exploring function spaces where the Ingham function is balanced,
either trivially or in a non-obvious way. Such an approach could facilitate
the transfer of information between these spaces, potentially paving
the way towards a proof of the Riemann hypothesis. The rich arithmetic
properties of the Ingham function revealed in this paper underscore
the viability of this program and suggest it as a compelling direction
for further exploration.


\begin{thebibliography}{10}
\bibitem[1]{key-1} Agarwal, R.P., Difference Equations and Inequalities,
Theory, methods, and applications, second ed., New York: Marcel Dekker,
Inc., (2000)

\bibitem[2]{key-2} Balazard, M., Le théorème des nombres premiers,
édition NANO, Calvage \& Mounet, (2016)

\bibitem[3]{key-3} Bender E. A. ; Goldman J. R., On the Applications
of Mobius Inversion in Combinatorial Analysis, The American Mathematical
Monthly, Vol. 82, No. 8 (1975)

\bibitem[4]{key-4} Bingham N.H. ; Inoue A., Tauberian and Mercerian
theorem for systems of Kernels, J. Math. Anal. Appl. 252, (2000)

\bibitem[5]{key-5} Bordellès, O., Arithmetic Tales, Advanced Edition,
Springer, (2019)

\bibitem[6]{key-6} Broughan, K., Equivalents of the Riemann hypothesis
(I), Arithmetic Equivalents, Cambridge University Press, (2017)

\bibitem[7]{key-7} Broughan, K., Equivalents of the Riemann hypothesis
(II), Analytic Equivalents, Cambridge University Press, (2017)

\bibitem[8]{key-8} Cloitre, B., Good variation theory: a Tauberian
approach to the Riemann hypothesis, International Journal of Mathematics
and Computer Science, Vol.2, (2016)

\bibitem[9]{key-9}Cohen H.; Dress F.; El Marraki M., Explicit estimates
for summatory functions linked to the Möbius, Funct. Approx. Comment.
Math. 37(1): 51-63 (2007)

\bibitem[10]{key-10} Diblík, J.; Schmeidel, E., On the existence
of solutions of linear Volterra difference equations, Appl. Math.
Comput., (2012)

\bibitem[11]{key-11} Erdös, P. ; Segal, S.L., A note on Ingham Summation
Method, Journal of Number Theory, (1978)

\bibitem[12]{key-12} Euler, L., Letter to Goldbach, 13 October 1729,
Euler Archive {[}E00715{]}, eulerarchive.maa.org, (1729)

\bibitem[13]{key-13} Hardy, G. H.; Littlewood, J. E. , Tauberian
theorems concerning power series and Dirichlet's series whose coefficients
are positive, Proceedings of the London Mathematical Society, (1914)

\bibitem[14]{key-14} Hardy, G.H.; Ramanujan, S., Twelve lectures
on subjects suggested by his life and work, AMS Chelsea Pub., pp 67-82,
(1999)

\bibitem[15]{key-15} Hindry, M., \href{https://www.cmls.polytechnique.fr/xups/xups12-02.pdf}{La preuve par André Weil de l'hypothèse de Riemann pour une courbe sur un corps fini},(2012)

\bibitem[16]{key-16} Ikehara, S., An extension of Landau's theorem
in the analytic theory of numbers, Journal of Mathematics and Physics
of the MIT, (1931)

\bibitem[17]{key-17}Ingham, A. E., Some Tauberian theorems connected
with the prime number theorem, J. London Math. Soc, 20 (1945)

\bibitem[18]{key-18} Jukes, K.A., On the Ingham and (D, h(n)) Summation
Methods, Journal of the London Mathematical Society, (1971)

\bibitem[19]{key-19} Kalmynin ,A.B.; Kosenko, P.R., Orthorecursive
expansion of unity, International Journal of Number Theory, Vol. 16,
(2020)

\bibitem[20]{key-20} Korevaar, J., Tauberian theory: A century of
developments, Springer, Volume 329, (2004)

\bibitem[21]{key-21} Levinson, N., On closure problems and the zeros
of the Riemann zeta function. PAMS 7, (1956)

\bibitem[22]{key-22}Montgomery, H., Vaughan, R.C., Multiplicative
number theory I, Cambridge University Press, (2006)

\bibitem[23]{key-23} Nyman B., On some groups and semigroups of translations,
Thesis, Uppsala university, (1950) 

\bibitem[24]{key-24} On-Line Encyclopedia of Integer Sequences, published
electronically at https://oeis.org

\bibitem[25]{key-25} Pennington, W. B., On Ingham summability and
summability by Lambert series, Proc. Camb. Philos. Soc., 51, (1955)

\bibitem[26]{key-26} J. Nelson Raja, P. Jaganathan and S. Domnic,
A New Variable-Length Integer Code for Integer Representation and
Its Application to Text Compression, Indian Journal of Science and
Technology, Vol 8 (24), (2015)

\bibitem[27]{key-27} Rajagopal, C. T., A note on Ingham summability
and summability by Lambert series, Proc. Indian Acad. Sci., A42, (1955)

\bibitem[28]{key-28} Segal, S.L., On Ingham's Summation Method, Canadian
Journal of Mathematics , Volume 18, (1966)

\bibitem[29]{key-29} Segal, S.L., Ingham\textquoteright s summability
method and Riemann hypothesis, Proc. London Math. Soc., (1975)

\bibitem[30]{key-30} Selberg, A., Old and new conjectures and results
about a class of Dirichlet series, In Proceedings of the Amalfi Conference
on Analytic Number Theory, (1992)

\bibitem[31]{key-31} Tenenbaum, G., Introduction à la théorie analytique
et probabiliste des nombres, collections SMF, (1995)

\bibitem[32]{key-32} Wiener, N., Tauberian theorems, Ann. of Math.
vol. 33. (1932)

\bibitem[33]{key-33} Wintner, A., Eratosthenian Averages, Baltimore,
Waverly Press, (1943)

\bibitem[34]{key-34} Zacharovas, V., A Tauberian theorem for the
Ingham summation method, Acta Arithmetica, (2011) 

\end{thebibliography}
\end{document}